\nonstopmode
\documentclass[11pt,draft,reqno]{amsart}
\input epsf.tex

\usepackage[all]{xy}
\usepackage{amsthm,array,amssymb,amscd,amsfonts,latexsym}
\xyoption{poly}

\DeclareFontFamily{OT1}{wncyr}{\hyphenchar\font45 }
\DeclareFontShape{OT1}{wncyr}{m}{n}{%
   <5> <6> <7> <8> <9> gen * wncyr
   <10> <10.95> <12> <14.4> <17.28> <20.74>  <24.88>wncyr10}{}
\DeclareFontShape{OT1}{wncyr}{m}{it}{%
   <5> <6> <7> <8> <9> gen * wncyi
   <10> <10.95> <12> <14.4> <17.28> <20.74> <24.88> wncyi10}{}
\DeclareFontShape{OT1}{wncyr}{m}{sc}{%
   <5> <6> <7> <8> <9> <10> <10.95> <12> <14.4>
   <17.28> <20.74> <24.88>wncysc10}{}
\DeclareFontShape{OT1}{wncyr}{b}{n}{%
   <5> <6> <7> <8> <9> gen * wncyb
   <10> <10.95> <12> <14.4> <17.28> <20.74> <24.88>wncyb10}{}
\input cyracc.def
\def\rus{\usefont{OT1}{wncyr}{m}{n}\cyracc\fontsize{9}{11pt}\selectfont}

\DeclareMathSizes{9}{9}{7}{5} 

\begin{document}

\theoremstyle{plain}
\newtheorem*{lemma*}{Lemma}
\newtheorem{lemma}[subsection]{Lemma}
\newtheorem*{theorem*}{Theorem}
\newtheorem{theorem}[subsection]{Theorem}
\newtheorem*{proposition*}{Proposition}
\newtheorem{proposition}[subsection]{Proposition}
\newtheorem*{corollary*}{Corollary}
\newtheorem{corollary}[subsection]{Corollary}
\theoremstyle{definition}
\newtheorem*{definition*}{Definition}
\newtheorem{definition}[subsection]{Definition}
\newtheorem*{example*}{Example}
\newtheorem{example}[subsection]{Example}
\theoremstyle{remark}
\newtheorem*{remark*}{Remark}
\newtheorem{remark}[subsection]{Remark}
\newtheorem{reduction}[subsection]{Reduction}
\newtheorem{problem}[subsection]{Problem}

\newcommand{\Span}{\operatorname{Span}}
\newcommand{\Symp}{\mbox{\boldmath$\rm Sp$}}
\newcommand{\g}{\mathfrak{g}}
\newcommand{\el}{\mathfrak{l}}
\newcommand{\lt}{\mathfrak{t}}
\newcommand{\lc}{\mathfrak{c}}
\newcommand{\lu}{\mathfrak{u}}
\newcommand{\lr}{\mathfrak{r}}
\newcommand{\Id}{\operatorname{id}}
\newcommand{\id}{\operatorname{id}}
\newcommand{\pr}{\operatorname{pr}}
\newcommand{\Hom}{\operatorname{Hom}}
\newcommand{\sign}{\operatorname{sign}}

\newcommand{\bbA}{{\mathbb A}}
\newcommand{\bbC}{{\mathbb C}}
\newcommand{\bbZ}{{\mathbb Z}}
\newcommand{\bbP}{{\mathbb P}}
\newcommand{\bbQ}{{\mathbb Q}}
\newcommand{\bbG}{{\mathbb G}}
\newcommand{\bbN}{{\mathbb N}}
\newcommand{\bbF}{{\mathbb F}}

\newcommand{\e}{{\lambda}}

\newcommand{\ve}{{\varepsilon}}
\newcommand{\vp}{{\varpi}}

\newcommand{\kbar}{\overline k}

\newcommand{\mo}{\mathopen<}
\newcommand{\mc}{\mathclose>}

\newcommand{\Ker}{\operatorname{Ker}}
\newcommand{\Aut}{\operatorname{Aut}}
\newcommand{\sdp}{\mathbin{{>}\!{\triangleleft}}} 
\newcommand{\Alt}{\operatorname{A}}   
\newcommand{\GL}{\mbox{\boldmath$\rm GL$}}

\newcommand{\PGL}{\mbox{\boldmath$\rm PGL$}}
\newcommand{\SL}{\mbox{\boldmath$\rm SL$}}
\newcommand{\rank}{\operatorname{rank}}
\newcommand{\aut}{\operatorname{Aut}}
\newcommand{\Char}{\operatorname{\rm char\,}} 
\newcommand{\diag}{\operatorname{\rm diag}}
\newcommand{\Gal}{\operatorname{Gal}}
\newcommand{\galois}{\Gal}
\newcommand{\lra}{\longrightarrow}
\newcommand{\SO}{\mbox{\boldmath$\rm SO$}}
\newcommand{\M}{\operatorname{M}}        
\newcommand{\ord}{\mathop{\rm ord}\nolimits}
\newcommand{\Sym}{{\operatorname{S}}}    
\newcommand{\tr}{\operatorname{\rm tr}}
\newcommand{\trace}{\tr}
\newcommand{\ad}{\operatorname{Ad}}

\newcommand{\Res}{\operatorname{Res}}
\newcommand{\Sha}{\mbox{\rus{\fontsize{11}{11pt}\selectfont{SH}}}}
\newcommand{\G}{\mathcal{G}}
\renewcommand{\H}{\mathcal{H}}
\newcommand{\gen}[1]{\langle{#1}\rangle}
\renewcommand{\O}{\mathcal{O}}
\newcommand{\C}{\mathcal{C}}
\newcommand{\Ind}{\operatorname{Ind}}
\newcommand{\End}{\operatorname{End}}
\newcommand{\Spin}{\mbox{\boldmath$\rm Spin$}}
\newcommand{\T}{\mathbf G}
\newcommand{\GT}{\mbox{\boldmath$\rm T$}}
\newcommand{\Inf}{\operatorname{Inf}}
\newcommand{\Tor}{\operatorname{Tor}}
\newcommand{\m}{\mbox{\boldmath$\mu$}}

\newcommand{\A}{{\sf A}}

\newcommand{\D}{{\sf D}}

\newcommand{\Lbd}{{\sf \Lambda}}


\title[On polarizations in
invariant theory]{On polarizations in invariant theory}

\author{Mark~Losik}
\address{M.~Losik: Saratov State University, ul. Astrakhanskaya, 83,
410026 Saratov, Russia} \email{LosikMV@info.sgu.ru}

\author{Peter~W.~Michor}
\address{P.~W.~Michor: Institut f\"ur Mathematik, Universit\"at Wien,
Strudlhofgasse 4, A-1090 Wien, Austria; {\it and:}
Erwin Schr\"odinger Institut f\"ur Mathematische
Physik, Boltzmanngasse 9, A-1090 Wien, Austria}
\email{Peter.Michor@esi.ac.at}

\author{Vladimir L. Popov}
\address{V.~L.~Popov: Steklov Mathematical Institute,
Russian Academy of Sciences, Gubkina 8, Moscow 119991,
Russia} \email{popovvl@orc.ru}

\thanks{M.\;L. and P.\;W.\;M. were supported
     by ``Fonds zur F\"orderung der
     wissenschaftlichen For\-schung,
     Projekt P~14195~MAT''. V.\;L.\;P. was
     supported in part by ETH, Z\"urich,
 Switzerland, Russian grants {\rus
N{SH}--123.2003.01} and {\rus RFFI 05--01--00455}, and
Program of Mathematics Section of Russian Academy of
Sciences.}
\date{May 4, 2005}


\begin{abstract}
Given a reductive algebraic group $G$ and a finite
dimensional algebraic $G$-module $V$, we study how
close is the algebra of $G$-invariant polynomials on
$V^{\oplus n}$  to the subalgebra generated by
polarizations of $G$-invariant polynomials on $V$. We
address this problem in a more general setting of
$G$-actions on arbitrary affine varieties.
\end{abstract}

\maketitle

\section {\bf Introduction}
\label{sect.intro}

\subsection{} Let $G$ be a reductive
algebraic group over an algebraically closed field $k$
of characteristic $0$, and let $V$ be a finite
dimensional algebraic $G$-module. Given a positive
integer $n$, consider the $G$-module $V^{\oplus
n}:=V\oplus\ldots\oplus V$ ($n$ summands). Finding
generators of the invariant algebra $k[V^{\oplus n}]^G$
of $V^{\oplus n}$ is the classical problem of invariant
theory. The classical method of constructing elements
of $k[V^{\oplus n}]^G$ is taking the polarizations
of invariants $f\in k[V]^G$, i.e., the polynomial
functions $f_{i_1,\ldots, i_n}$ on $V^{\oplus n}$ given
by the formal expansions
\begin{equation}\textstyle \label{pola}
f(x_1v_1+\ldots+x_nv_n)=\sum_{i_1,\ldots,i_n\in \mathbb
Z_+} x_1^{i_1}\ldots x_n^{i_n}
f_{i_1,\ldots,i_n}(v_1,\ldots, v_n),
\end{equation}
where $(v_1,\ldots, v_n)$ is generic element of
$V^{\oplus n}$ and $x_1,\ldots, x_n$ are variables. Let
${\rm pol}_n k[V]^G$ be the subalgebra of $k[V^{\oplus
n}]^G$ generated by the polarizations of all the $f$'s.

There are $G$-modules enjoying the property
\begin{equation}\label{equal}
{\rm pol}_n k[V]^G=k[V^{\oplus n}]^G.\end{equation} For
instance, \eqref{equal} holds, by Study's theorem,
\cite{study}, for  the standard action
 of  $G={\bf O}_m$ on $V=k^m$.
 By Weyl's theorem, \cite{weyl},
\eqref{equal} holds for $G={\rm S}_m$ acting on $V=k^m$
by permuting the coordinates. In \cite{hunziker},
\eqref{equal} is established for the natural action of
the Weyl group $G$ of type ${\sf B}_m$ on $V=k^m$ and
for the standard action of the dihedral group $G$ on
$V=k^2$.

However, in general, ${\rm pol}_n k[V]^G$ and
$k[V^{\oplus n}]^G$ do not coincide. For instance, for
the natural action of $G={\bf SL}_n$ on $V=k^n$ clearly
${\rm pol}_n k[V]^G=k$  (since $k[V]^G=k$) but
$k[V^{\oplus n}]^G\neq k$. It is less easy to find
examples where \eqref{equal} fails for finite $G$, but
such examples exist as well: in \cite{wallach} it was
observed that \eqref{equal} does not hold for the
natural action of the Weyl group $G$ of type ${\sf
D}_m$ on $V=k^m$ ($m\geqslant 4$) for $n\geqslant 2$.

In this paper we analyze the relationship between
$k[V^{\oplus n}]^G$ and ${\rm pol}_n k[V]^G$. We prove
that if $G$ is finite, then $k[V^{\oplus n}]^G$ is the
integral closure of ${\rm pol}_n k[V]^G$ in its field
of fractions, and the natural morphism of affine
varieties determined by these algebras is bijective.
Actually, instead of linear actions we consider the
more general setting of actions on arbitrary affine
varieties for which we define a generalization of
polarizations. In this setting, we prove that if $G$ is
finite, then
  the invariant algebra is integral over the
subalgebra generated by  generalized polarizations, and
the natural dominant morphism between affine varieties
determined by these algebras is injective (in the
graded case, bijective).

For connected $G$, one cannot expect such results, as
the example of ${\bf SL}_n$ acting on $k^n$ shows. This
naturally leads to distinguishing the $n$'s for which
$k[V^{\oplus n}]^G$ is integral over ${\rm pol}_n
k[V]^G$ and defining {\it the polarization index of}
$V$,
$${\rm pol\hskip .7mm ind}\hskip .4mm(V),$$
as the supremum taken over all such $n$'s. We prove
that $k[V^{\oplus m}]^G$ is integral over ${\rm pol}_m
k[V]^G$ for every $m\leqslant {\rm pol\hskip .7mm
ind}\hskip .4mm(V)$, and show that calculating ${\rm
pol\hskip .7mm ind}\hskip .4mm(V)$ is closely related
to the old problem of describing linear subspaces lying
in the Hilbert nullcone of~$V$ (see \cite{hilbert},
\cite{mufoki}, \cite{popov-vinberg},
 \cite{gerstenhaber}, \cite{crt},
\cite{mor}, \cite{mufoki}, and the references therein),
namely, to analyzing a certain geometric property of
such subspaces.

Using this reduction, we calculate the polarization
index of some $G$-modules $V$. Namely, we prove that if
$G$ is a finite group or a linear algebraic torus, then
${\rm pol\hskip .7mm ind}\hskip .4mm(V)=\infty$. For
$G={\bf SL}_2$, we describe all linear subspaces of $V$
lying in the Hilbert nullcone of $V$ and prove that
${\rm pol\hskip .7mm ind}\hskip .4mm(V)=\infty$ if $V$
does not contain a simple $2$-dimensional submodule,
and ${\rm pol\hskip .7mm ind}\hskip .4mm(V)=1$
otherwise. Finally, we calculate the polarization index
of every semisimple Lie algebra $\mathfrak g$: we prove
that ${\rm pol\hskip .7mm ind}\hskip .4mm(\mathfrak
g)=1$ if $\mathfrak g$ is not isomorphic to $\mathfrak
{sl}_2\oplus\ldots\oplus \mathfrak {sl}_2$, and ${\rm
pol\hskip .7mm ind}\hskip .4mm(\mathfrak g)=\infty$
otherwise. As an application to the above mentioned old
topic of linear subspaces lying in the Hilbert
nullcone, we prove that a semisimple Lie algebra
$\mathfrak g$ contains a $2$-dimensional nilpotent
nontriangularizable linear subspace if and only if
$\mathfrak g$ is not isomorphic to $\mathfrak
{sl}_2\oplus\ldots\oplus \mathfrak {sl}_2$.

\subsection{Notation}

$k[X]$ is the algebra of regular functions of an
algebraic variety $X$. If $X$ is irreducible,  $k(X)$
is the field of rational function of $X$.

If a group $S$ acts on a set $Z$, we put $Z^S:=\{z\in
Z\mid s\cdot z=z\ \mbox{for all}\ s\in S\}$.

Below every action of an algebraic group is algebraic
(morphic).

$G^0$ is the identity component of an algebraic group
$G$.

If $X$ is an affine variety endowed with an action of a
reductive algebraic group $G$, then $\pi^{}_{X, G}:
X\to X/\!\!/G$ is the categorical quotient, i.e.,
$X/\!\!/G$ is an affine algebraic variety and
$\pi^{}_{X, G}$ a dominant (actually, surjective)
morphism such that $\pi^{*}_{X,
G}(k[X/\!\!/G])=k[X]^G$.

Given  a linear algebraic torus $T$, its character
group ${\rm Hom}(T, {\bf G}_m)$ is written additively.
The value of $\lambda\in {\rm Hom}(T, {\bf G}_m)$ at
$t\in T$ is denoted by $t^\lambda$.  For an algebraic
$T$-module $V$ (not necessarily finite dimensional),
$V_\lambda$ is the $\lambda$-isotypic component of $V$,
\begin{equation*} V_\lambda:=\{v\in V\mid t\cdot v=
t^\lambda v \quad\mbox{for every $t\in T$}\}.
\end{equation*}

By $ \langle v_1,\ldots, v_n\rangle $ we denote the
linear span of vectors $v_1,\ldots, v_n$ of a vector
space over $k$.

We set $\mathbb N:=\{1, 2, \ldots\}$ and $\mathbb
Z_+:=\{0, 1, 2, \ldots\}$.

\section{\bf Generalized polarizations}

\subsection{}
Let a reductive algebraic group $G$ act on the
irreducible affine algebraic varieties $X$ and $Y$. Let
$Z$ be an irreducible affine algebraic variety endowed
with an action of a linear algebraic torus $T$. The set
$\Lambda:=\{\lambda\in {\rm Hom}(T, {\bf G}_m)\mid
k[Z]_\lambda\neq 0\}$ is then a submonoid of ${\rm
Hom}(T, {\bf G}_m)$ and the isotypic components yield a
$\Lambda$-grading of $k[Z]$:
\begin{equation}\textstyle\label{oplus}
k[Z]=\bigoplus_{\lambda\in \Lambda}k[Z]_\lambda,\quad
k[Z]_\mu k[Z]_\nu\subseteq k[Z]_{\mu +\nu}\ \mbox{for
all $\mu, \nu\in \Lambda$}.
\end{equation}

The groups $G$ and $T$ act on $Y\times Z$ through the
first and second factors, respectively. From $k[Y\times
Z]=k[Y]\otimes k[Z]$ and \eqref{oplus} we obtain
\begin{equation}\textstyle\label{decomp1}
k[Y\times Z]=\bigoplus_{\lambda\in \Lambda} k[Y]\otimes
k[Z]_\lambda\quad\mbox{and}\quad k[Y\times
Z]^G=\bigoplus_{\lambda\in \Lambda} k[Y]^G\otimes
k[Z]_\lambda.
\end{equation}
We identify $k[Y]$ and $k[Z]$ respectively with the
subalgebras $k[Y]\otimes 1$ and $1\otimes k[Z]$ of
$k[Y\times Z]$.

 Assume now
that there is an open $T$-orbit $\mathcal O$ in $Z$.
This condition is equivalent to either of the following
properties ${\rm (o1)}$, ${\rm (o2)}$, see
\cite[Theorem 3.3]{popov-vinberg}:
\begin{gather}
{{\rm (o1)}}\hskip 3mm k(Z)^T=k;\notag\\
\label{multiplicity} {{\rm (o2)}}\hskip 2mm\dim
k[Z]_\lambda= 1\hskip 2mm \mbox{for every $\lambda\in
\Lambda$}.
\end{gather}
For every $\lambda\in \Lambda$, fix a nonzero element
$b_\lambda\in k[Z]_\lambda$. Multiplying every
$b_\lambda$ by an appropriate scalar we may assume that
\begin{equation}\label{monoid}
b_\mu b_\nu=b_{\mu+\nu}\quad\mbox{for all $\mu, \nu\in
\Lambda$}.
\end{equation}
Indeed, fix a point $x_0\in \mathcal O$. The definition
of $k[Z]_\lambda$ implies that
 $b_\lambda(x_0)\neq 0$, so replacing $b_\lambda$ by
$b_\lambda/b_\lambda(x_0)$ we may assume that
$b_\lambda(x_0)=1$. Then \eqref{monoid} follows
from~\eqref{multiplicity},~\eqref{oplus}.

 From \eqref{decomp1} and
\eqref{multiplicity} it follows that every  $h\in
k[Y\times Z]$ admits a unique decomposition
\begin{equation}\textstyle\label{decomp2}
h=\sum_{\lambda\in\Lambda} p_{\lambda} b_{\lambda},
\quad p_{\lambda}\in k[Y]
\end{equation}
(in \eqref{decomp2} all but finitely many $p_\lambda$'s
are equal to zero), and $h$ lies in $k[Y\times Z]^G$ if
and only if $p_\lambda\in k[Y]^G$ for all $\lambda$.
From \eqref{monoid} we obtain
\begin{equation}\textstyle\label{product}
\bigl(\sum_{\mu\in \Lambda}p'_\mu b_\mu\bigl)
\bigl(\sum_{\nu\in \Lambda}p''_\nu b_\nu\bigl)=
\sum_{\lambda\in \Lambda}\bigl(\sum_{\mu+\nu=\lambda}
p'_\mu {p''}_\nu\bigr) b_\lambda, \quad\mbox{$p'_\mu,
p''_\nu\in k[Y]$}.
\end{equation}

Consider now a $G$-equivariant morphism
\begin{equation}\label{phi}
\varphi: Y\times Z\longrightarrow X.
\end{equation}
\begin{definition}\label{generalized} Let $f\in k[X]^G$. The invariants $p_\lambda\in
k[Y]^G$ defined   by {\rm \eqref{decomp2}} for
$h=\varphi^*(f)$ are called {\it the
$\varphi$-polarizations of $f$}. The subalgebra of
$k[Y]^G$ generated by all  the $\varphi$-polarizations
of the elements of $k[X]^G$ is denoted by ${\rm
pol}_{\varphi} k[X]^{G}$ and called {\it the
$\varphi$-polarization algebra of~$Y$}.
\end{definition}

\begin{remark} More generally, if $\Phi$ is a
collection of $G$-equivariant morphisms \eqref{phi}
(where $Z$ and $X$ depend on $\varphi$), then one can
define {\it the $\Phi$-polarization algebra  of $Y$} as
the subalgebra  of $k[Y]^G$ generated by all the
$\varphi$-polarization algebras of $Y$ for
$\varphi\in\Phi$. \end{remark}

Since changing the $b_\lambda$'s clearly replaces the
$\varphi$-polarizations of $f\in k[X]^G$ by their
scalar multiples, the algebra ${\rm pol}_{\varphi}
k[X]^{G}$ does not depend on the choice of the
$b_\lambda$'s.

\begin{example}\label{point} If $Z$ is a single
point, \eqref{phi} is a morphism $\varphi: Y\rightarrow
X$, and ${\rm pol}_{\varphi}
k[X]^{G}=\varphi^*(k[X]^G)$.
\end{example}
\begin{example} \label{classical-setting1}(Classical setting) Let $V$ be a finite dimensional
algebraic $G$-module and let $n\in \mathbb N$. Take
$X=V$ and $Y=V^{\oplus n}$ with the diagonal
$G$-action. Let $Z$ be ${\bf A}^n$ endowed with the
natural action of the diagonal torus $T$ of ${\bf
GL}_n$,
\begin{equation*}
{\rm diag}(t_1,\ldots, t_n)\cdot (\alpha_1,\ldots,
\alpha_n)= (t_1\alpha_1,\ldots, t_n\alpha_n).
\end{equation*}
\end{example}
Identifying $(i_1,\ldots,  i_n)\in \mathbb Z^n$ with
the character $T\rightarrow {\bf G}_m$, ${\rm
diag}(t_1,\ldots, t_n)\mapsto t_1^{i_1}\cdots
t_n^{i_n}$, we identify $\mathbb Z^n$ with ${\rm
Hom}(T, {\bf G}_m)$. Then $\Lambda=\mathbb Z_+^n$. If
$z_1,\ldots, z_n$ are the standard coordinate functions
on $Z$, then for every $\lambda=(i_1,\ldots, i_n)\in
\Lambda$, the isotypic component $k[Z]_{\lambda}$ is
spanned by $b_\lambda:=z_1^{i_1}\cdots z_n^{i_n}$. So,
condition \eqref{multiplicity} holds. Clearly,
\eqref{monoid} holds as well.

 Recall that the classical
$n$-polarizations of a polynomial $f\in k[V]$ are the
polynomials $f_{i_1,\ldots, i_n}\in k[V^{\oplus n}]$,
where $(i_1,\ldots, i_n)\in \Lambda$, such that
\begin{equation}\textstyle\label{calssical}
f(\sum_{j=1}^n\alpha_jv_j)=\sum_{i_1,\ldots, i_n\in
\Lambda}\alpha_1^{i_1}\cdots \alpha_n^{i_n}
f_{i_1,\ldots, i_n}(v_1,\ldots, v_n)\quad\mbox{for all
$v_j\in V$, $\alpha_j\in k$}.
\end{equation}

Since $\alpha_1^{i_1}\cdots \alpha_n^{i_n}$ is the
value of $z_1^{i_1}\cdots z_n^{i_n}$ at
$(\alpha_1,\ldots, \alpha_n)\in Z$,  it readily follows
from \eqref{calssical} and Definition~\ref{generalized}
that the classical $n$-polarizations of $f$ are the
$\varphi$-polarizations of $f$ for
\begin{equation}\label{tau}
\varphi:=\tau_n: Y\times Z\rightarrow
X,\quad\bigl((v_1,\ldots, v_n), (\alpha_1,\ldots,
\alpha_n)\bigr)\mapsto\alpha_1v_1+\ldots+ \alpha_nv_n.
\end{equation}
In this setting, we denote the $\varphi$-polarization
algebra of $Y$ by ${\rm pol}_nk[V]^{G}$. \quad
$\square$

\begin{example}\label{ex3} If $G={\bf O}_m$ and $V=k^m$ with the natural $G$-action, then
${\rm pol}_nk[V]^G=k[V^{\oplus n}]^G$ by Study's
theorem,~\cite{study}.

If $G={\bf Sp}_m$, $m$ even, and $V=k^m$ with the
natural $G$-action, then $ {\rm pol}_{\tau_n\times
\tau_n}k[V^{\oplus 2}]^G=k[V^{\oplus n}]^G$ (see
\eqref{tau}) by \cite{weyl}.

If $G={\bf SL}_m$ and $V=k^m$ with the natural
$G$-action, then ${\rm pol}_{\tau_n\times\ldots\times
\tau_n} k[V^{\oplus m}]^G= k[V^{\oplus n}]^G$,
see\;\cite{weyl}. \quad $\square$
\end{example}

\vskip 1.5mm

 From \eqref{product} we deduce that the algebra ${\rm
pol}_{\varphi} k[X]^{G}$ is generated by all
$\varphi$-polarizations of the $f$'s for $f$ running
through the generators of $k[X]^G$. Since by Hilbert's
theorem, the algebra $k[X]^G$ is finitely generated,
this means that the algebra ${\rm pol}_{\varphi}
k[X]^{G}$ is finitely generated as well. Hence there is
an affine algebraic variety that we denote by
$Y/\!\!/\varphi$, and a dominant morphism $\pi_\varphi:
Y\to Y/\!\!/\varphi$ such that
$\pi_\varphi^*(k[Y/\!\!/\varphi])= {\rm pol}_{\varphi}
k[X]^{G}$. Since ${\rm pol}_{\varphi} k[X]^{G}\subseteq
k[Y]^G$, the definition of categorical quotient for the
$G$-action on $Y$ implies that there is a dominant
morphism $\nu: Y/\!\!/G\to Y/\!\!/\varphi$ such that
the following diagram is commutative:
\begin{equation}\label{diagram}
\begin{matrix}\xymatrix @=4.2mm@M=
1.5mm@R=7mm{Y\ar[rr]^{\pi^{}_{Y,
G}}\ar[rd]_{\pi_\varphi}&&Y/\!\!/G\ar[ld]^\nu
\\
&Y/\!\!/\varphi & }
\end{matrix}\quad .
\end{equation}

The set of all morphisms from $Z$ to $X$ is endowed
with the $G$-action defined by the formula $(g\cdot
\psi)(z):=g\cdot(\psi(z))$ for $\psi:Z\to X$, $g\in G$,
$z\in Z$. Using \eqref{phi}, we can consider $Y$ as a
$G$-stable algebraic family of such morphisms. Namely,
with every $y\in Y$ we associate the morphism
\begin{equation}\label{phiy}
\varphi_y: Z\longrightarrow X,\hskip 2mm z\mapsto
\varphi(y,z).
\end{equation}
Then for every $z\in Z$ and $g\in G$ we have
$\varphi_{g\cdot y}(z)=\varphi(g\cdot y, z)=
\varphi(g\cdot(y, z))=g(\varphi(y, z))= (g\cdot
\varphi_y) (z)$, so $\varphi_{g\cdot y}= g\cdot
\varphi_y$.
\begin{lemma}\label{fibers} For every $y_1, y_2\in Y$,
the following properties are equivalent:
\begin{enumerate}
\item[\rm (i)] $\pi_\varphi^{-1}(\pi_\varphi(y_1))=
\pi_\varphi^{-1}(\pi_\varphi(y_2))$; \item[\rm (ii)]
$\pi^{}_{X, G}\circ \varphi_{y^{}_1}= \pi^{}_{X,
G}\circ \varphi_{y^{}_2}$.
\end{enumerate}
\end{lemma}
\begin{remark}\label{separation} Property (i)
means that points $y_1, y_2\in Y$ are not separated by
the $\varphi$-polarization algebra ${\rm pol}_{\varphi}
k[X]^{G}$.
\end{remark}
\begin{proof} By virtue of \eqref{phiy}, property
(ii) is equivalent to the property  \begin{equation}
\label{points} \pi^{}_{X, G}(\varphi(y_1,z))=
\pi^{}_{X, G}(\varphi(y_2,z))\quad\mbox{for all $z\in
Z$}.
\end{equation}
Since the variety $X$ is affine, for a fixed $z\in Z$,
equality in \eqref{points} holds if and only if
\begin{equation}\label{f}
s(\pi^{}_{X, G}(\varphi(y_1,z)))= s(\pi^{}_{X,
G}(\varphi(y_2,z)))\quad \mbox{for every $s\in
k[X/\!\!/G]$}.
\end{equation}
Since $\pi_{X, G}^*(k[X/\!\!/G])=k[X]^G$, in turn,
\eqref{f}
 is equivalent to the property
\begin{equation}\label{ff}
f(\varphi(y_1,z))= f(\varphi(y_2,z))\quad \mbox{for
every $f\in k[X]^G$}.
\end{equation}
Setting $h=\varphi^*(f)$ for $f$ in \eqref{ff}, we thus
obtain $h(y_1, z)=h(y_2, z)$ for all $z\in Z$, i.e.,
using the notation of \eqref{decomp2},
$\sum_{\lambda\in
\Lambda}p_\lambda(y_1)b_\lambda=\sum_{\lambda\in
\Lambda}p_\lambda(y_2)b_\lambda$. Since $\{b_\lambda\}$
are linearly independent, this shows that the equality
in \eqref{ff} is equivalent to the collection of
equalities $p_\lambda(y_1)= p_\lambda(y_2)$,
$\lambda\in \Lambda$. Definition~\ref{generalized} and
Remark~\ref{separation} now imply the claim. \quad
$\square$
\renewcommand{\qed}{}
\end{proof}

\begin{lemma}\label{orbits} If $G$ is a finite group, then
for every two morphisms $\psi_i: Z\rightarrow X$, $i=1,
2$, the following properties are equivalent:
\begin{enumerate}
\item[\rm(i)] $\pi^{}_{X, G}\circ \psi_1=\pi^{}_{X,
G}\circ \psi_2$; \item[\rm(ii)] there is $g\in G$ such
that $\psi_2=g\cdot\psi_1$.
\end{enumerate}
\end{lemma}
\begin{proof} (ii)$\Rightarrow$(i) is clear
(and holds for every reductive $G$, not necessarily
finite). Assume now that (i) holds. Consider in
$Z\times X$ the closed subset
\begin{equation}\label{Psi}
\Psi:=\{(z, x)\in Z\times X \mid \pi^{}_{X,
G}(\psi_1(z))=\pi^{}_{X, G}(x)\}.
\end{equation}

Since $G$ is finite, every fiber of $\pi^{}_{X, G}$ is
a $G$-orbit, see, e.g.,\;\cite[Theorem
4.10]{popov-vinberg}. Hence for $(z, x)\in Z\times X$,
the condition $\pi^{}_{X, G}(\psi_1(z))=\pi^{}_{X,
G}(x)$ in \eqref{Psi} is equivalent to the existence of
$g\in G$ such that $x
=g\cdot(\psi_1(z))=(g\cdot\psi_1)(z)$. In turn, the
last equality means that the point $(z,x)$ lies in the
graph of  $g\cdot\psi_1$,
\begin{equation}\label{graph}
\Gamma_{g\cdot \psi_1}:=\{\bigl(z, (g\cdot
\psi_1)(z)\bigr)\in Z\times X\mid z\in Z\}.
\end{equation}
 On the other hand,
\eqref{Psi}, \eqref{graph} clearly imply that
$\Gamma_{g\cdot \psi_1}\subseteq \Psi$ for every $g$.
Thus,
\begin{equation}\textstyle\label{cup}
\Psi=\bigcup_{g\in G}\Gamma_{g\cdot\psi_1}.
\end{equation}
But every $\Gamma_{g\cdot\psi_1}$ is a closed subset of
$Z\times X$ isomorphic to $Z$. So, by \eqref{cup},
$\Psi$ is a union of finitely many closed irreducible
subsets of the same dimension. Hence these subsets are
precisely the irreducible components of $\Psi$.

On the other hand, it follows from (i) that
\begin{equation}\label{Psi2}
\Psi:=\{(z, x)\in Z\times X \mid \pi^{}_{X,
G}(\psi_2(z))=\pi^{}_{X, G}(x)\}.
\end{equation}
Using the above argument, we then deduce from
\eqref{Psi2} that the graph of $\psi_2$,
\begin{equation*}
\Gamma_{\psi_2}:=\{(z, \psi_2(z))\in Z\times X\mid z\in
Z\},
\end{equation*}
 is an irreducible component of
$\Psi$ as well. Therefore there is $g\in G$ such that
$\Gamma_{\psi_2}=\Gamma_{g\cdot\psi_1}$. Hence $g\cdot
\psi_1=\psi_2$, i.e., (ii) holds. \quad $\square$
\renewcommand{\qed}{}
\end{proof}

\begin{theorem} \label{finite1} Maintain the
notation of this section. If $G$ is a finite group,
then \begin{enumerate} \item[\rm(i)] the morphism $\nu:
Y/\!\!/G\rightarrow Y/\!\!/\varphi$ in \eqref{diagram}
is injective; \item[\rm(ii)]  $k(Y)^G$ is the field of
fractions of the $\varphi$-polarization algebra ${\rm
pol}_{\varphi} k[X]^{G}$.
\end{enumerate}
\end{theorem}

\begin{proof} Since $G$ is finite, fibers of
$\pi^{}_{X, G}$ are precisely $G$-orbits. On the other
hand, by Lem\-mas~\ref{fibers} and \ref{orbits}, every
fiber of $\pi_\varphi$ is a $G$-orbit as well. This and
the commutative diagram \eqref{diagram} yield (i).
Since $\nu$ is dominant and ${\rm char}\;k=0$, from (i)
it follows that $\nu$ is a birational isomorphism.
Since $Y/\!\!/\varphi$ is affine, $k(Y/\!\!/\varphi)$
is the field of fractions of $k[Y/\!\!/\varphi]$. This
and \eqref{diagram} now imply (ii). \quad $\square$
\renewcommand{\qed}{}
\end{proof}

\subsection{} Under a supplementary assumption there is a
geometric criterion of finiteness of $\nu$. It is based
on a general statement essentially due to Hilbert.
Namely, consider an action of a reductive algebraic
group $G$ on an irreducible affine algebraic variety
$M$. Assume that the corresponding $G$-action on $k[M]$
preserves a $\mathbb Z_+$-grading $k[M]=\oplus_{n\in
\mathbb Z_+}k[M]_n$ such that $k[M]_0=k$ and $\dim
k[M]_n<\infty$ for every $n$. Let $A$ be a homogeneous
subalgebra of $k[M]^G$.
\begin{lemma}\label{hilbert} The following
properties are equivalent:
\begin{enumerate}
\item[\rm(i)] $\{x\in M\mid f(x)=0\hskip 1.2mm
\forall\hskip .5mm f\in \oplus_{n\in \mathbb
N}k[M]^G_n\}=\{x\in M\mid h(x)=0\hskip 1.2mm
\forall\hskip .5mm h\in \oplus_{n\in \mathbb N}A_n\}$;
\item[\rm (ii)] $k[M]^G$ is integral over $A$.
\end{enumerate}
If these properties hold and $G$ is connected, then
$k[M]^G$ is the integral closure of $A$ in $k[M]$.
\end{lemma}
\begin{proof}
For linear actions, (i)$\Rightarrow$(ii) is proved by
Hilbert in \cite[\S4]{hilbert}. In the general case the
argument  is the same. Implication (ii)$\Rightarrow$(i)
is clear. The last statement follows from the first
since it is well known that $k[M]^G$ is integrally
closed in $k[M]$ for connected $G$ (connectedness of
$G$ implies that $G$ acts trivially on the set of roots
of the equation of integral dependence). \quad
$\square$
\renewcommand{\qed}{}
\end{proof}
Lemma~\ref{hilbert} implies the following geometric
criterion of finiteness of $\nu$.
 Assume that the
 $G$-actions on $X$ and $Y$
 can be
 extended
 to the $G\times {\bf G}_m$-actions
 such
 that
\begin{gather}\label{equivar}
\mbox{$\varphi$ is
 $G\times {\bf G}_m$-equivariant,}\\
\label{k} k[X]^{{\bf G}_m}=k[Y]^{{\bf G}_m}=k.
\end{gather}
 From \eqref{k} we then deduce that  ${\rm Hom}({\bf
G}_m, {\bf G}_m)$ can be identified with $\mathbb Z$ so
that the isotypic component decompositions of $k[X]$
and $k[Y]$ become the $\mathbb Z_+$-gradings of these
algebras,
\begin{equation}\textstyle\label{grading}
k[X]=\bigoplus_{n\in \mathbb Z_+} k[X]_n, \quad
k[Y]=\bigoplus_{n\in \mathbb Z_+} k[Y]_n.
\end{equation}
Since  every isotypic component is a finitely generated
module over invariants,  see, e.g.,\;\cite[Theorem
3.24]{popov-vinberg}, from \eqref{k} we deduce that
these gradings enjoy the properties
\begin{equation}
k[X]_0=k[Y]_0=k \hskip 3mm \mbox{and}\hskip 2mm \dim
k[X]_n<\infty,\ \dim k[Y]_n<\infty\ \mbox{for every
$n\in \mathbb Z_+$}.
\end{equation}

It follows from \eqref{equivar} that $k[X]^G$ and
$k[Y]^G$ are graded subalgebras of the graded algebras
$k[X]$ and $k[Y]$ respectively, and from
Definition~\ref{generalized} we deduce that ${\rm
pol}_\varphi k[X]^G$ is a graded subalgebra of the
graded algebra $k[Y]^G$.

 The ideal $\oplus_{n\in \mathbb N} k[X]_n$  in $k[X]$
(respectively, $\oplus_{n\in \mathbb N} k[Y]_n$ in
$k[Y]$) is maximal and ${\bf G}_m$-stable, so the point
$0_X\in X$ (respectively, $0_Y\in Y$) where it
vanishes, is ${\bf G}_m$-fixed. As invariants separate
closed orbits, see, e.g.,\,\cite[Theorem
4.7]{popov-vinberg}, \eqref{k} implies that $X^{{\bf
G}_m}=\{0_X\}$, $Y^{{\bf G}_m}=\{0_Y\}$. Hence $0_X\in
X^G$, $0_Y\in Y^G$. From \eqref{equivar} we deduce that
$\varphi(0_Y\times Z)=0_X$. We put
\begin{equation}\label{nullcone}
\mathcal N^{}_{Y, G}:=\pi^{-1}_{Y, G}(\pi^{}_{Y,
G}(0_Y)),\quad \mathcal P^{}_{Y, G}:=\pi^{-1}_{\varphi}
(\pi^{}_{\varphi}(0_Y)),\quad \mathcal N^{}_{X,
G}:=\pi^{-1}_{X, G}(\pi^{}_{X, G}(0_X)).
\end{equation}
By virtue of \eqref{diagram}, the following inclusion
holds:
\begin{equation}\label{subseteq}
\mathcal N^{}_{Y, G} \subseteq \mathcal P^{}_{Y, G}.
\end{equation}
Since $\mathcal N^{}_{Y, G}$ is precisely the set of
points of $Y$ whose $G$-orbit contains $0_Y$ in the
closure,
\begin{equation}\label{S}
\mathcal N^{}_{S, G}=S\cap \mathcal N^{}_{Y, G}
\end{equation}
for every $G$-stable closed subset $S$ of $Y$
containing $0_Y$.

\begin{example} \label{classical-setting2} Maintain the
notation of Example~\ref{classical-setting1}. Then the
${\bf G}_m$-actions on $X=V$ and $Y=V^{\oplus n}$ by
scalar multiplications yield the $G\times {\bf
G}_m$-extensions of $G$-actions such that
\eqref{equivar}, \eqref{k} hold. Thus, in the classical
setting, the assumptions of Subsection 2.2 hold. In
this case, $0_X=0$, $0_Y=(0,\ldots, 0)$. The varieties
$\mathcal N^{}_{V^{\oplus n}, G}$ and $\mathcal
N^{}_{V, G}$ are respectively the Hilbert nullcones of
$G$-modules $V^{\oplus n}$ and $V$, and $\mathcal
P^{}_{V^{\oplus n}, G}$ is the locus of the maximal
homogeneous ideal of ${\rm pol}_nk[V]^G$.\quad
$\square$
\renewcommand{\qed}{}
\end{example}

\begin{lemma}\label{graded} Maintain the assumptions
of Subsection {\rm 2.2}. The following  properties are
equiva\-lent:
\begin{enumerate}
\item[\rm(i)] $\nu$ is finite; \item[\rm(ii)] $\mathcal
N^{}_{Y, G}= \mathcal P^{}_{Y, G}$.
\end{enumerate}
\end{lemma}
\begin{proof} This immediately follows from
Lemma~\ref{hilbert}. \quad $\square$
\renewcommand{\qed}{}
\end{proof}

\begin{theorem} \label{finite-clas}
Maintain the assumptions of Subsection {\rm 2.2} and
let $G$ be finite. Then
\begin{enumerate}
\item[\rm(i)] $\nu$ is finite and bijective;
\item[\rm(ii)] if $Y$ is normal, $\nu:
Y/\!\!/G\rightarrow Y/\!\!/\varphi$ is the
normalization of $Y/\!\!/\varphi$, and $k[Y]^G$ is the
integral closure of ${\rm pol}_{\varphi} k[X]^{G}$ in
$k(Y)^G$; \item[\rm(iii)] if $Y$ is normal, ${\rm
pol}_{\varphi} k[X]^{G}=k[Y]^G$ if and only if ${\rm
pol}_{\varphi} k[X]^{G}$ is integrally closed.
\end{enumerate}
\end{theorem}
\begin{proof}
Theorem~\ref{finite1}(i) implies that $\mathcal
N^{}_{Y, G}=\mathcal P^{}_{Y, G}$ ($=0_Y$). Hence $\nu$
is finite by Lemma~\ref{graded}. Being finite, $\nu$ is
closed, and since $\nu$ is also dominant,
Theorem~\ref{finite1}(i) implies that $\nu$ is
bijective. This proves (i). If $Y$ is normal, then
$Y/\!\!/G$ is normal as well, see, e.g.,\;\cite[Theorem
3.16]{popov-vinberg}. Since by
Theorem~\ref{finite1}(ii), $\nu$ is a birational
isomorphism, this, (i), and the definitions of
$Y/\!\!/G$, $Y/\!\!/\varphi$, $\nu$ prove (ii). Claim
(iii) follows from (ii).
 \quad $\square$
\renewcommand{\qed}{}
\end{proof}
\begin{corollary} In the classical setting $($see
Example~{\rm\ref{classical-setting1}}$)$, let $G$ be
finite. Then \begin{enumerate} \item[\rm(i)] $\nu$ is
finite and bijective; \item[\rm (ii)] $k[V^{\oplus
n}]^G$ is the integral closure of ${\rm pol}_n
k[V]^{G}$ in $k(V^{\oplus n})^G$; \item[\rm(iii)]
$k[V^{\oplus n}]^G={\rm pol}_n k[V]^{G}$ if and only if
${\rm pol}_n k[V]^G$ is integrally closed.
\end{enumerate}
\end{corollary}
\begin{example} Maintain the notation of
Example~{\rm \ref{classical-setting1}} and let $V=k^m$.
If $G$ is the symmetric group in $m$ letters acting on
$V$ by permuting the coordinates, then $k[V^{\oplus
n}]^G={\rm pol}_nk[V]^{G} $ for every $n$, \cite{weyl}.
This equality also holds for the Weyl group of type
${\sf B}_m$ and the dihedral groups, \cite{hunziker}.
But for the Weyl group of type ${\sf D}_m$, $m\geqslant
4$, and $n=2$ it does not hold,~\cite{wallach}.

Namely, ${\sf D}_m$ acts on the standard coordinate
functions $x_1,\ldots, x_m$ on $V$ by permutations and
changes of an even number of signs, and $k[V]^{{\sf
D}_m}=k[\sigma_1,\ldots, \sigma_m]$ where
\begin{equation*}\textstyle
\sigma_s=\sum_{i=1}^m x_i^{2s}\hskip3mm \mbox{for
$1\leqslant s\leqslant m-1$},\hskip5mm
\sigma_m=x_1\cdots x_m,
\end{equation*}
see, e.g.,\;\cite{humphreys}. Take another copy of $V$
with the standard coordinate functions $y_1,\ldots,
y_m$, and naturally identify $x_1,\ldots, x_m,
y_1,\ldots, y_m$ with the functions on $V^{\oplus 2}$.
Then $k[V^{\oplus 2}]^{{\sf D}_m}$ is generated by
${\rm pol}_2k[V]^{{\sf D}_m}$ and the polynomials
\begin{equation}\textstyle\label{P}
P_{r_1}\dots P_{r_d}(\sigma_n),\quad r_1,\ldots, r_d
\hskip2mm \mbox{odd},
\quad \sum_{i=1}^dr_i\leqslant n-d,
\end{equation}
where $P_r:=\sum_{i=1}^m y_i^r\frac{\partial}{\partial
x_i}$, see \cite{wallach}, \cite{hunziker}. The group
${\sf B}_m$ is generated by ${\sf D}_m$ and the
reflection $w$ such that $w\cdot x_i=x_i$ for $i<m$ and
$w\cdot x_m=-x_m$. The operators $P_{r_i}$ from
\eqref{P} commute with the diagonal action of ${\sf
B}_m$ on $V^{\oplus 2}$, therefore $w(P_{r_1}\dots
P_{r_d}(\sigma_n))=-P_{r_1}\dots P_{r_d}(\sigma_n)$.
This yields
\begin{equation}\label{equation}
(P_{r_1}\dots P_{r_d}(\sigma_n))^2\in k[V^{\oplus
2}]^{{\sf B}_m}. \end{equation} Since $k[V^{\oplus
2}]^{{\sf B}_m}={\rm pol}_2k[V]^{{\sf B}_m}$ and,
clearly, ${\rm pol}_2k[V]^{{\sf B}_m}\subseteq {\rm
pol}_2k[V]^{{\sf D}_m}$, we deduce from
\eqref{equation} that $k[V^{\oplus 2}]^{{\sf D}_m}$ is
integral over ${\rm pol}_2k[V]^{{\sf D}_m}$. This
agrees with Theorem~\ref{finite-clas} (that gives more
delicate information).
\end{example}

\section{\bf Polarization index}

\subsection{} In this section we take up the classical
setting and maintain the notation of
Examples~\ref{classical-setting1},
\ref{classical-setting2}, and that of \eqref{nullcone}.
If $n, m\in \Bbb N$ and $n\leqslant m$, we naturally
identify $V^{\oplus n}$ with the subspace
$\{(v_1,\ldots, v_n,0,\ldots,0)\mid v_i\in V\}$ of
$V^{\oplus m}$. It is then not difficult to see that
\begin{equation}\label{P}
\mathcal P_{V^{\oplus n}, G}=V^{\oplus n}\cap \mathcal
P_{V^{\oplus m}, G}.
\end{equation}

\begin{lemma}\label{subspace} The following
properties of  a point $v=(v_1,\ldots, v_n)\in
V^{\oplus n}$ are equivalent: \begin{enumerate}
\item[\rm(i)] $v\in \mathcal P^{}_{V^{\oplus n}\!, G}$;
\item[\rm(ii)] $\langle v_1,\ldots, v_n\rangle
\subseteq \mathcal N^{}_{V, G}$.
\end{enumerate}
\end{lemma}
\begin{proof}
Let $f\in k[V]^G$ be a nonconstant homogeneous
function. If $v \in \mathcal P^{}_{V^{\oplus n}\!, G}$,
then the definition of $\mathcal P^{}_{V^{\oplus n}\!,
G}$ (see  \eqref{nullcone} and
Example~\ref{classical-setting2}) yields that, in the
notation of \eqref{pola}, we have $f_{i_1,\ldots,
i_n}(v)=0$ for all $i_1,\ldots,i_n\in \mathbb Z_+$.
 From this and \eqref{pola}
we obtain
\begin{equation}\label{sum}
f(\alpha_1v_1+\ldots+\alpha_nv_n)=0\quad \mbox{\rm for
all $\alpha_i\in k$}.
\end{equation}
So,  $\langle v_1,\ldots, v_n\rangle$ lies in the zero
set of every $f$. The definition of $\mathcal N^{}_{V,
G}$ (see  \eqref{nullcone} and
Example~\ref{classical-setting2}) now implies that
$\langle v_1,\ldots, v_n\rangle \subseteq \mathcal
N^{}_{V\!, G}$.

Conversely, assume that the last inclusion holds. By
the definition of $\mathcal N^{}_{V, G}$, this implies
\eqref{sum}. By \eqref{pola}, this in turn yields that
$f_{i_1,\ldots, i_n}(v)=0$ for all $i_1,\ldots,i_n\in
\mathbb Z_+$. The definition of $\mathcal
P^{}_{V^{\oplus n}, G}$ then implies that $v \in
\mathcal P^{}_{V^{\oplus n}\!, G}$.
 \quad $\square$
\renewcommand{\qed}{}
\end{proof}

\begin{definition}\label{index} {\it The polarization
index of a $G$-module} $V$ is
\begin{equation*}\label{pol}
{\rm pol}\;{\rm ind}\;(V):={\rm sup}\;n
\end{equation*}
with the supremum taken over all $n$ such that in
\eqref{subseteq} the equality holds, $\mathcal
N^{}_{V^{\oplus n}\!, G}=\mathcal P^{}_{V^{\oplus n}\!,
G}$.
\end{definition}

From \eqref{nullcone}, Definition~\ref{index}, and the
equality ${\rm pol}_1 k[V]^G=k[V]^G$ we obtain
\begin{equation}\label{>0}
{\rm pol{\hskip .8mm}ind}\;(V)\geqslant 1.
\end{equation}
It is also clear that
\begin{equation}\label{trivial}
{\rm pol\;ind}(V\oplus U)={\rm pol\;ind}(V)\hskip
2mm\mbox{if $U$ is a trivial $G$-module}.
\end{equation}

\begin{lemma} \label{np} For every $n\in \mathbb N$,
\begin{equation*}
\mathcal N^{}_{V^{\oplus n}\!, G}
\begin{cases}=\mathcal P^{}_{V^{\oplus n}\!,
G}& \mbox{if \hskip 1mm $n\leqslant {\rm pol}\;{\rm ind}\;(V)$},\\
\varsubsetneq \mathcal P^{}_{V^{\oplus n}\!,G}&
\mbox{if \hskip 1mm $n> {\rm pol}\;{\rm ind}\;(V)$}.
 \end{cases}
\end{equation*}
\end{lemma}

\begin{proof} By virtue of Definition~\ref{index},
for $n> {\rm pol}\;{\rm ind}\;(V)$, this follows from
\eqref{subseteq}, and for $n\leqslant {\rm pol}\;{\rm
ind}\;(V)$, from \eqref{S} and \eqref{P}. \quad
$\square$
\renewcommand{\qed}{}
\end{proof}

\begin{corollary}\label{integral}
The extension ${\rm pol}_nk[V]^G\!\subseteq\!
k[V^{\oplus n}]^G$ is integral if and only if
$n\!\leqslant\! {\rm pol}\;{\rm ind}\hskip .2mm(V)$.
\end{corollary}
\begin{proof}
This follows from Lemmas~\ref{np} and \ref{graded}.
\quad $\square$
\renewcommand{\qed}{}
\end{proof}

Call a character ${\bf G}_m\rightarrow {\bf G}_m$,
$t\mapsto t^d$, {\it positive} if $d>0$. Every
homomorphism $\gamma: {\bf G}_m\rightarrow G$ endows
$V$ with the structure of ${\bf G}_m$-module defined by
$t\cdot v:=\gamma(t)\cdot v$. We denote by $V(\gamma)$
the submodule of this ${\bf G}_m$-module equal to the
sum of all the isotypic components whose weight is
positive. Clearly, if $v\in V(\gamma)$, then the
closure of ${\bf G}_m$-orbit (and, all the more,
$G$-orbit) of $v$ contains $0_V$. Hence
$V(\gamma)\subseteq \mathcal N^{}_{V, G}$. The
Hilbert--Mumford theorem, \cite{hilbert}, \cite{mufoki}
(see, e.g.,\;\cite[5.3]{popov-vinberg}), claims that
\begin{equation}\textstyle \label{hm}
\mathcal N^{}_{V, G}=\bigcup_{\gamma} V(\gamma).
\end{equation}

\begin{lemma}\label{weights}
The following properties of an integer $n\in \mathbb N$
are equivalent:
\begin{enumerate}
 \item[\rm (i)] for every linear subspace $L$
such that $\dim L\leqslant n$ and $L\subseteq \mathcal
N^{}_{V\!, G}$, there is a homomorphism $\gamma: {\rm
G}_m\rightarrow G$ such that $L\subseteq V(\gamma)$;
\item[\rm (ii)] $n\leqslant {\rm pol}\;{\rm ind}\;(V)$.
\end{enumerate}
\end{lemma}
\begin{proof} Let (i) hold. Take a point
$v=(v_1,\ldots, v_n)\in \mathcal P^{}_{V^{\oplus n}\!,
G}$. By Lemma~\ref{subspace},  $\langle v_1,\ldots,
v_n\rangle $ is contained in $\mathcal N^{}_{V\!, G}$.
By (i), $\langle v_1,\ldots, v_n\rangle \subseteq
V(\gamma)$ for some $\gamma$. This implies that the
closure of ${\bf G}_m$-orbit (and, all the more,
$G$-orbit) of $v$ contains $0^{}_{V^{\oplus n}}$, i.e.,
$v\in \mathcal N^{}_{V^{\oplus n}\!, G}$. So, by
\eqref{subseteq}, we have $\mathcal P^{}_{V^{\oplus
n}\!, G}= \mathcal N^{}_{V^{\oplus n}\!, G}$, whence
$n\leqslant {\rm pol}\;{\rm ind}\;(V)$ by
Definition~\ref{index}. This proves
(i)$\Rightarrow$(ii).

Conversely, let (ii) holds. Consider in $\mathcal
N^{}_{V\!, G}$ a linear subspace $L$ of dimension
$\leqslant n$. Then $L=\langle v_1,\ldots, v_n\rangle$
for some $v_i\in V$. By Lemma~\ref{subspace}, the point
$v=(v_1,\ldots, v_n)$ lies in  $\mathcal
P^{}_{V^{\oplus n}\!, G}$. By (ii), and Lemma~\ref{np},
we have $v\in \mathcal N^{}_{V^{\oplus n}\!, G}$. From
\eqref{hm} we now deduce that $v\in V^{\oplus
n}(\gamma)$ for some $\gamma$. Since $V^{\oplus
n}(\gamma)=V(\gamma)^{\oplus n}$, this yields $v_i\in
V(\gamma)$ for every $i$, or, equivalently, $L\subseteq
V(\gamma)$. Thus (ii)$\Rightarrow$(i) is proved. \quad
$\square$
\renewcommand{\qed}{}
\end{proof}

\begin{corollary} \label{UV} {\rm (i)}
Let $U$ be a submodule
of $V$. Then
\begin{equation*}{\rm pol{\hskip .9mm}ind}\;(U)\geqslant
{\rm pol{\hskip .9mm}ind}\;(V).
\end{equation*}
\begin{enumerate}
\item[\rm (ii)] Let $V_i$ be a finite dimensional
algebraic module of a reductive algebraic group $G_i$,
$i=1,\ldots, m$, and let
$G=G_1\times\ldots\times G_m$, 
$V=V_1\oplus\ldots\oplus V_m$. Then
\begin{equation*}\label{polii}
{\rm pol\;ind}(V)=\min_i {\rm pol\hskip .8mm ind}(V_i).
\end{equation*}
\end{enumerate}
\end{corollary}
\begin{proof}
Statement (i) readily follows from Lemma~\ref{weights}.
Let the assumptions of (ii) hold, and let
$\pi_i:V\rightarrow V_i$, $p_i: G\rightarrow G_i$ be
the natural projections. Since the Hilbert nullcone is
the set of points whose orbits contain zero in the
closure, we have $\mathcal N^{}_{V, G}=\mathcal
N^{}_{V_1, G_1}\times\ldots \times \mathcal N^{}_{V_m,
G_m}$. A linear subspace $L$ lying in $\mathcal
N^{}_{G, V}$ is contained in $V(\gamma)$ for some
$\gamma$ if and only if $\pi_i(L) \subseteq
V_i(p_i\circ\gamma)$ for every $i$. Using these
properties and Lemma~\ref{weights}, we deduce (ii).
\quad $\square$
\renewcommand{\qed}{}
\end{proof}

\subsection{} We now calculate the polarization index of
some $G$-modules.

\begin{theorem}\label{finite2} If $G$ is a finite
group, then for any $G$-module $V$,
\begin{equation*}
{\rm pol}\;{\rm ind}\;(V)=\infty.
\end{equation*}
\end{theorem}
\begin{proof} This follows from
Theorem~\ref{finite-clas}, Lemma~\ref{graded}, and
Definition~\ref{index}.\quad $\square$
\renewcommand{\qed}{}
\end{proof}

\begin{theorem}\label{Torus} If $G$ is
a linear algebraic torus, then for any $G$-module $V$,
\begin{equation*}
{\rm pol}\;{\rm ind}\;(V)=\infty.
\end{equation*}
\end{theorem}
\begin{proof} It is well known
(and immediately follows from \eqref{hm}) that in this
case there are homomorphisms $\gamma_i: {\bf
G}_m\rightarrow G$, $i=1,\ldots, s$, such that
\begin{equation}\label{torus}\mathcal N^{}_{V\!,
G}=V(\gamma_1)\cup\ldots \cup V(\gamma_s).
\end{equation}
Since every linear subspace $L$ of $V$ is an
irreducible algebraic variety, \eqref{torus} implies
that if $L\subseteq\mathcal N^{}_{V\!, G}$, then
$L\subseteq V(\gamma_i)$ for some $i$, whence the claim
by Lemma~\ref{weights}.\quad $\square$
\renewcommand{\qed}{}
\end{proof}

\begin{lemma}\label{trivinv} Let $k[V]^G$=k.
\begin{enumerate}
\item[\rm(i)] ${\rm pol}\;{\rm ind}\;(V)={\rm sup}\;n$
where the supremum is taken over all n such that
$k[V^{\oplus n}]^G=k$. \item[\rm (ii)] If $G^0$ is
semisimple, then ${\rm pol}\;{\rm ind}\;(V)$ is equal
to the generic transitivity degree of the $G$-action on
$V$, see \cite{popov}, i.e., to the maximum $n$ such
that there is an open $G$-orbit in $V^{\oplus n}$. In
this case,
$${\rm pol}\;{\rm ind}\;(V)\leqslant \dim G/\dim V.$$
\end{enumerate}
\end{lemma}
\begin{proof}
The condition $k[V]^G=k$ and
Definition~\ref{generalized} imply that  $\mathcal
P^{}_{V^{\oplus n}\, G}=V^{\oplus n}$ for every~$n$. On
the other hand, $k[V^{\oplus n}]^G=k$ is equivalent to
$\mathcal N^{}_{V^{\oplus n}\, G}=V^{\oplus n}$. This
gives (i). Being semisimple, $G^0$ has no nontrivial
characters, hence $k[V^{\oplus n}]^G=k$ is equivalent
to the existence of an open $G$-orbit in $V^{\oplus
n}$, see \cite[Theorem 3.3 and the Corollary of Theorem
2.3] {popov-vinberg}. This proves (ii).\quad $\square$
\renewcommand{\qed}{}
\end{proof}
\begin{example}
If $G={\bf SL}_m$ and $V=k^m$ with the natural
$G$-action, then Lemma~\ref{trivinv} implies ${\rm
pol}\;{\rm ind}\;(V)= m-1$.

If $G={\bf Sp}_{m}$ and $V=k^{m}$ ($m$ even) with the
natural $G$-action, then Lemma~\ref{trivinv} implies
${\rm pol}\;{\rm ind}\;(V)=1$.

If $G={\bf O}_m$ and $V=k^m$ with the natural
$G$-action, then Example~\ref{ex3} and
Definition~\ref{index} yield ${\rm pol}\;{\rm
ind}\;(V)=\infty$.

If $G={\bf SO}_m$ and $V=k^m$ with the natural
$G$-action, then the classical description of
$k[V^{\oplus n}]^G$, see \cite{weyl}, implies that
$k[V^{\oplus n}]^G$ is integral over $k[V^{\oplus
n}]^{{\bf O}_m}$. Hence in this case again ${\rm
pol}\;{\rm ind}\;(V)=\infty$, however, in contrast to
the case of ${\bf O}_m$,  the algebras $k[V^{\oplus
n}]^G$ and ${\rm pol}_nk[V]^G$ do not coincide if $m$
divides $n$.\quad $\square$
\end{example}

We now calculate the polarization index
 of any ${\bf SL}_2$-module.
Denote by $R_d$ the ${\bf SL}_2$-module of binary forms
in $x$ and $y$ of degree $d$, see, e.g.,
\cite[0.12]{popov-vinberg}. Up to isomorphism, $R_d$ is
the unique simple ${\bf SL}_2$-module of dimension
$d+1$. According to the classical Hilbert theorem,
\cite[\S5]{hilbert} (see, e.g.,\,\cite[Example 1 in
5.4]{popov-vinberg}),
\begin{equation} \textstyle \label{hilbert1}
\mathcal N^{}_{R_d\!, {\bf SL}_2}=\bigcup_{l\in
R_1}l^{[d/2]+1}R_{d-[d/2]-1},
\end{equation}
and for every nonzero $l\in R_1$, there is a
homomorphism $\gamma: {\bf G}_m\rightarrow {\bf SL}_2$
such that
\begin{equation}\label{hilbert2}
l^{[d/2]+1}R_{d-[d/2]-1}=R_d(\gamma)
\end{equation}
and vice versa.

\begin{lemma}\label{sl2subspaces} For $d\geqslant 2$,
the following properties of a linear subspace $L$ of
$R_d$ lying in $\mathcal N^{}_{R_d\!, {\bf SL}_2}$ are
equivalent:
\begin{enumerate}
\item[\rm(i)] $L$ is maximal $($\hskip -.1mm with
respect to inclusion\hskip .2mm$)$ among the linear
subspaces lying in $\mathcal N^{}_{R_d\!, {\bf SL}_2}$;
\item[\rm(ii)] there is $l\in R_1$, $l\neq 0$ such that
$L=l^{[d/2]+1}R_{d-[d/2]-1}$.
\end{enumerate}
\end{lemma}

\begin{proof} Using that
$k[x, y]$ is a unique factorization domain and every
$l\in R_1$, $l\neq 0$ is a simple element in it, we
obtain that for every nonzero $l_1, l_2\in R_1$,
\begin{equation}\label{0}
l_1^{[d/2]+1}R_{d-[d/2]-1}\cap
l_2^{[d/2]+1}R_{d-[d/2]-1}=\{0\}\quad \mbox{if
$l_1/l_2\notin k$}.
\end{equation}
Therefore
 it suffices to show that for every
 $2$-dimensional
 linear subspace $P$ lying in
 $\mathcal N^{}_{R_d\!, {\bf SL}_2}$ there is
 $l\in R_1$ such that $P\subseteq
l^{[d/2]+1}R_{d-[d/2]-1}$. Let $f_1, f_2$ be a basis of
$P$.
 Then \eqref{hilbert1} implies that
\begin{equation}\label{hilbert3}
f_i=l_i^{[d/2]+1}h_i\quad \mbox{for some $l_i\in R_1,
h_i\in R_{d-[d/2]-1}$}.
\end{equation}
We have to show that if $\alpha f_1+\beta f_2\in
\mathcal N^{}_{R_d\!, {\bf SL}_2}$ for every $\alpha,
\beta\in k$, then $l_1/l_2\in k$.

For contradiction, assume that $l_1$ and $l_2$ are
linearly independent. Applying ${\bf SL}_2$, we then
may assume that $l_1=x$, $l_2=y$. Since $P\subseteq
\mathcal N^{}_{R_d\!, {\bf SL}_2}$,  from
\eqref{hilbert1} we deduce that for every $\alpha,
\beta\in k$ there are $\mu, \nu\in k$, $h\in
R_{d-[d/2]-1}$ (depending on $\alpha, \beta$) such that
$ \alpha x^{[d/2]+1}h_1+\beta y^{[d/2]+1}h_2= (\mu
x+\nu y )^{[d/2]+1}h. $ Using that $k[x, y]$ is unique
factorization domain, we deduce from this equality and
\eqref{hilbert3} that $\mu\nu\neq 0$ if
 $\alpha\beta\neq 0$. Hence we may
  assume that for every nonzero $\alpha, \beta\in k$
  there are $\mu\in k$, $h\in
R_{d-[d/2]-1}$ (depending on $\alpha, \beta$) such that
\begin{equation}\label{relation1}
\alpha x^{[d/2]+1}h_1+\beta y^{[d/2]+1}h_2= (\mu x+y
)^{[d/2]+1}h.
\end{equation}

Note that
\begin{equation}\label{range} \mbox{when
$\alpha$ and $\beta$ in \eqref{relation1} vary, $\mu$
ranges over an infinite set}.
\end{equation}
Indeed, otherwise \eqref{relation1} implies that there
is a basis of $P$ whose elements are divisible by some
$(\mu x +y)^{[d/2]+1}$. Hence $x^{[d/2]+1}h_1$ is
divisible by $(\mu x +y)^{[d/2]+1}$ as well. Since
$\deg h_1<[d/2]+1$, this is impossible.

We now consider separately the cases of even and odd
$d$. First, let $d$ be even, $d=2m$. Then
\begin{equation} \textstyle \label{h1}
h_1=\sum_{i=0}^{m-1}\eta_ix^{m-i-1}y^i,\quad \eta_i\in
k.
\end{equation}
Plugging \eqref{h1} in equality \eqref{relation1}, and
then differentiating it $m$ times with respect to $x$,
substituting  $y=-\mu x$, and dividing both sides by
$\alpha x^m$, we obtain the following equality:
\begin{equation} \textstyle \label{equation1}
\sum_{i=0}^{m-1}(-1)^i\frac{(2m-i)!}{(m-i)!}\;\eta_i
\mu^i=0.
\end{equation}
 Since
$h_1\neq 0$, \eqref{equation1} contradicts
\eqref{range}.

Let now $d$ be odd, $d=2m-1$. Then $h_1$ is still given
by \eqref{h1} and
\begin{equation}\textstyle \label{h2}
h_2=\sum_{j=0}^{m-1}\theta_jx^{j}y^{m-j-1},\quad
\theta_j\in k.
\end{equation}
Plugging \eqref{h1}, \eqref{h2} in equality
\eqref{relation1}, and then differentiating it $m-1$
times respectively with respect to $x$ and $y$,
substituting $y=-\mu x$, and dividing both sides by
$x^m$, we obtain respectively the equalities
\begin{gather}\textstyle \label{1}
\theta_{m-1} \mu^m=\frac{ \alpha}{ \beta}
\sum_{i=0}^{m-1}(-1)^{m+i-1}\frac{(2m-i-1)!}
{(m-1)!(m-i)!}\;\eta_i
\mu^i,\\
\textstyle\label{2} \alpha
\eta_{m-1}+\beta\sum_{j=0}^{m-1}(-1)^{m-j}
\frac{(2m-j-1)!}{(m-1)!(m-j)!}\;\theta_j \mu^{m-j}=0.
\end{gather}
Multiplying \eqref{2} by $\theta_{m-1}\mu^{m-1}$,
replacing $\theta_{m-1}\mu^m$ by the right hand side of
\eqref{1}, and dividing both sides by $\alpha$, we
obtain
\begin{equation}\textstyle\label{3}
\begin{gathered}\textstyle
\eta_{m-1}\theta_{m-1} \mu^{m-1} + \sum_{i,
j=0}^{m-1}(-1)^{i-j-1}
\frac{(2m-j-1)!(2m-i-1)!}{((m-1)!)^2(m-i)!(m-j)!}\;\eta_i
\theta_j \mu^{m-j+i-1}=0.
\end{gathered}
\end{equation}
 From  \eqref{range} we deduce that all the coefficients
of the left hand side of \eqref{3}, considered as a
polynomial in $\mu$, vanish. In particular,
\begin{equation}\label{coeff}
\eta_0\theta_{m-1}=\eta_{m-1}\theta_0=0.
\end{equation}

If $\theta_{m-1}=0$, then \eqref{1}, \eqref{range}
imply $\eta_0=\ldots=\eta_{m-1}=0$ contrary to $h_1\neq
0$. Similarly, if $\eta_{m-1}=0$, then \eqref{2}
\eqref{range} imply $\theta_0=\ldots=\theta_{m-1}=0$
contrary to $h_2\neq 0$. Thus,
$\eta_{m-1}\theta_{m-1}\neq 0$, whence, by
\eqref{coeff}, $\eta_0=\theta_0=0$. From \eqref{h1},
\eqref{h2}, \eqref{relation1} we then deduce that for
$m\geqslant 2$, the left hand side of \eqref{relation1}
is divisible by $xy$. Hence $h$ in \eqref{relation1} is
divisible by $xy$ as well; in particular, $m\geqslant
3$. Thus, for $m\geqslant 3$, dividing both sides of
\eqref{relation1} by $xy$, we obtain $\alpha x^{m-1}
h'_1+\beta y^{m-1} h'_2= (\mu x+y)^{m-1} h', $ with
$h'_1, h'_2, h\in R_{m-2}$. This means that in
considering \eqref{relation1} we may step down from
case $m$ to case $m-1$. Continuing this way we reduce
the consideration of \eqref{relation1} to the case
$m=2$. In this case, the above argument shows that $h$
is a nonzero element of $R_1$ divisible by $xy$. This
contradiction completes the proof.
 \quad $\square$
\renewcommand{\qed}{}
\end{proof}

\begin{corollary}
The action of ${\bf SL}_2$ on the set of maximal
linear subspaces of $R_d$ lying in $\mathcal N^{}_{R_d,
{\bf SL}_2}$ is transitive. The dimension of every such
subspace is equal to $d-[d/2]$ and $2$ respectively for
for $d\neq1$ and $d=1$.
\end{corollary}

\begin{theorem}\label{sl2} For $G={\bf SL}_2$ and
$V=R_{d_1}\oplus\ldots\oplus R_{d_m}$,
\begin{equation*}
{\rm pol}\;{\rm ind}\;(V)=\begin{cases}1&\mbox{if
$d_i=1$ for some $i$},\\
\infty&\mbox{otherwise}.\end{cases}
\end{equation*}
\end{theorem}
\begin{proof}  Since the $G$-module
$R_0$ is trivial, by \eqref{trivial} we may assume that
$d_i\geqslant 1$ for every $i$. Since $k[R_1]^{{\bf
SL}_2}=k$ and $k[R_1^{\oplus 2}]^{{\bf SL}_2}\neq k$,
Lemma~\ref{trivinv} implies that ${\rm pol\hskip .7mm
ind}\hskip .4mm(R_1)=1$.
 From this, the Corollary
of Lemma~\ref{weights}, and \eqref{>0} we deduce the
claim for the cases where $d_i=1$ for some~$i$.

Assume now that $d_i\geqslant 2$ for every $i$, and let
$L$ be a linear subspace of $V$ lying in $\mathcal
N^{}_{V, G}$. Let $\pi_i: V=R_{d_1}\oplus\ldots\oplus
R_{d_m}\rightarrow R_{d_i}$ be the natural projection
to the $i$th summand, and let $L_i:=\pi_i(L)$. Since
$\pi_i(\mathcal N^{}_{V, G})= \mathcal N^{}_{R_i, G}$,
we have $L_i\subset \mathcal N^{}_{R_i, G}$ for every
$i$. Hence by Lemma~\ref{sl2subspaces} and
\eqref{hilbert2}, for every $i$, there is homomorphism
$\gamma_i:{\bf G}_m\rightarrow G$ such that
$L_i\subseteq R_{d_i}(\gamma_i)$. Take now a point
$v\in L$ such that $\pi_i(v)\neq 0$ for every $i$.
Since $L$ lies in $\mathcal N^{}_{V, G}$, it follows
from \eqref{hm} that $v\in V(\gamma)$ for some
$\gamma$. Hence $\pi_i(v)\in R_{d_i}(\gamma)$. But
$\pi_i(v)\in R_{d_i}(\gamma_i)$ as well. By
\eqref{hilbert2}, \eqref{0}, this yields
$R_{d_i}(\gamma_i)=R_{d_i}(\gamma)$. Hence $L\subseteq
V(\gamma)$. From Lemma~\ref{weights} we now deduce that
${\rm pol\hskip .7mm ind}\hskip .4mm(V)= \infty$
completing the proof.
  \quad $\square$
\renewcommand{\qed}{}
\end{proof}

\begin{theorem} Let $G$ be a connected semisimple
algebraic group and let $\mathfrak g$ be its Lie
algebra endowed with the adjoint $G$-action. Then
\begin{equation*}
{\rm pol\hskip .7mm ind}\hskip .4mm(\mathfrak g)=
\begin{cases}1
&\mbox{if\hskip 3mm $\mathfrak g$ is not isomorphic to
$\mathfrak{sl}_2\oplus\ldots\oplus
\mathfrak{sl}_2$},\\
\infty&\mbox{ otherwise}.
\end{cases}
\end{equation*}
\end{theorem}
\begin{proof} In this case,
$\mathcal N^{}_{\mathfrak g, G}$ is the cone of all
nilpotent elements in $\mathfrak g$,
see,\;e.g.,\;\cite[5.1]{popov-vinberg}. Every subspace
$\mathfrak g(\gamma)$ is the unipotent radical of a
parabolic subalgebra of $\mathfrak g$, see\;\cite[VIII,
4.4]{bourbaki2}, \cite[8.4.5]{springer}, and hence lies
in a maximal (with respect to inclusion) unipotent
subalgebra of $\mathfrak g$. Maximal unipotent
subalgebras of $\mathfrak g$ are precisely the
unipotent radicals of Borel subalgebras of $\mathfrak
g$, and $G$ acts transitively on the set of such
subalgebras, see, e.g.,\;\cite[Ch.\;6]{springer}. This
implies that for
 a linear
subspace $L$  of $\mathfrak g$ lying in $\mathcal
N^{}_{\mathfrak g, G}$ the following properties are
equivalent:
\begin{enumerate}
\item[(i)] the subalgebra of $\mathfrak g$ generated by
$L$ is unipotent (i.e., lies in $\mathcal
N^{}_{\mathfrak g, G}$); \item[(ii)] there is a
homomorphism $\gamma: {\bf G}_m\rightarrow G$ such that
$L\subseteq \mathfrak g(\gamma)$.
\end{enumerate}
 From this, \eqref{>0},  and Lemma~\ref{weights} we
deduce that equality ${\rm pol\hskip .7mm ind}\hskip
.4mm\mathfrak g=1$ is equivalent to the following
property:
 there is a $2$-dimensional
linear subspace $L$ of $\mathfrak g$ such that $L$ lies
in $\mathcal N^{}_{\mathfrak g, G}$ but the subalgebra
of $\mathfrak g$ generated by $L$ does not lie in
$\mathcal N^{}_{\mathfrak g, G}$. If this property
holds, we say, for brevity, that $\mathfrak g$ is a
$2$-algebra.

We shall show now that if
 $\mathfrak g$ is not isomorphic to
 $\mathfrak{sl}_2\oplus\ldots\oplus
\mathfrak{sl}_2$, then $\mathfrak g$ is a $2$-algebra.
To this end we remark that if a semisimple subalgebra
$\mathfrak s$ of $\mathfrak g$ is a $2$-algebra, then
$\mathfrak g$ is a $2$-algebra as well: since the cone
of nilpotent elements of $\mathfrak s$ lies in
$\mathcal N^{}_{\mathfrak g, G}$, this readily follows
from the definition of a $2$-algebra. Given this
remark, we see that the following two statements
immediately imply our claim:
\begin{enumerate}
\item[(a)] if $\mathfrak g\ncong
\mathfrak{sl}_2\oplus\ldots\oplus\mathfrak {sl}_2$,
then $\mathfrak g$ contains a subalgebra isomorphic to
either $\mathfrak {sl}_3$ or $\mathfrak {so}_5$;
\item[(b)] $\mathfrak {sl}_3$ and $\mathfrak {so}_5$
are $2$-algebras
\end{enumerate}
(note that in $\mathfrak {so}_5$ there are no
subalgebras isomorphic to $\mathfrak {sl}_3$, and vice
versa).

To prove (a), denote by $\Phi$ the root system of
$\mathfrak g$ with respect to a fixed maximal torus.
Let $\alpha_1,\ldots, \alpha_l$ be a system of simple
roots in $\Phi$ (enumerated as in \cite{bourbaki1}).
Fix a Chevalley basis $\{X_\alpha, X_{-\alpha},
H_\alpha\}_{\alpha\in \Phi}$ of $\mathfrak g$,
\cite{bourbaki2}. We may assume that $\mathfrak g$ is
simple, $\mathfrak g\ncong \mathfrak {sl}_2$,
$\mathfrak {so}_5$. For such $\mathfrak g$,
  it is easily seen that
  there are two roots $\lambda, \mu\in
\Phi$ such that the subalgebra of $\mathfrak g$
generated by $X_\lambda$ and $X_\mu$ is isomorphic to
$\mathfrak{sl}_3$: for $\mathfrak g$ of types ${\sf
A}_l$ $(l\geqslant 2)$, ${\sf B}_l$ $(l\geqslant 3)$,
${\sf C}_l$ $(l\geqslant 3)$, ${\sf D}_l$ $(l\geqslant
4)$, ${\sf F}_4$, one can take $\lambda=\alpha_1$,
$\mu=\alpha_2$; for types ${\sf E}_6$, ${\sf E}_7$,
${\sf E}_8$, take $\lambda=\alpha_1$, $\mu=\alpha_3$;
for type ${\sf G}_2$, take $\lambda=\alpha_2$, $\mu=
3\alpha_1+\alpha_2$. This proves (a).

We turn now to the proof of (b). In $\mathfrak {sl}_3$
we explicitly present a subspace $L$ enjoying the
desired properties (we are grateful to H.~Radjavi for
this example, \cite{ra}). Namely, take
\begin{gather}\label{sl3}
L:=\langle X_{\alpha_1}\!+\!X_{\alpha_2},
X_{-\alpha_1}-X_{-\alpha_2}\rangle=
\Bigl\{\!\begin{bmatrix}0&a&0\\
b&0&a\\
0&-b&0\end{bmatrix}\mid a, b\in k\Bigr\}.
\end{gather}
Then  \eqref{sl3} implies that the subalgebra generated
by $L$ contains the element
$H_{\alpha_1}-H_{\alpha_2}$. Since it  is semisimple,
this subalgebra does not lie in $\mathcal
N^{}_{\mathfrak g, G}$. On the other hand, the matrix
in the right hand side of \eqref{sl3} is nilpotent
(this is equivalent to the property that the sums of
all its principal minors of orders $2$ and $3$ are
equal to $0$, and this is immediately verified). So,
 $L\subseteq\mathcal N^{}_{\mathfrak
g, G}$. This proves (b) for $\mathfrak {sl}_3$.

Let now $\mathfrak g=\mathfrak {so}_5$. In this case,
an explicit construction of the desired subspace $L$ is
unknown to us, so we shall use an indirect argument.
The underlying space of $\mathfrak g$ is the space of
all skew-symmetric $5\times5$-matrices. Let $x_{ij}\in
k[\mathfrak g]$, $1\leqslant i,j\leqslant 5$, be the
standard coordinate functions on $\mathfrak g$ given by
$x_{ij}\bigl((a_{pq})\bigr)=a_{ij}$. Then
$x_{ij}=-x_{ji}$. Consider the matrix $A:=(x_{ij})$.
Then $k[\mathfrak g]^G=k[f_2, f_4]$ where $f_2$, $f_4$
are the coefficients of the characteristic polynomial
of $A$, i.e., ${\rm det}(tI_5-A)= t^5+f_2t^3+f_4t$,
see, e.g.,\;\cite{popov-vinberg}. The Newton formulas
expressing the sums of squares of eigenvalues of $A$
via the elementary symmetric functions of them imply
that ${\rm tr} (A^2)=-2f_2$, ${\rm
tr}(A^4)=2f_2^2-4f_4$. Hence
\begin{equation}\label{trace}
k[\mathfrak g]^G=k[{\rm tr} (A^2), {\rm tr}(A^4)].
\end{equation}

Let now $y_{ij}, z_{ij}\in k[{\mathfrak g}^{\oplus
2}]$, $1\leqslant i,j\leqslant 5$, be the standard
coordinate functions on ${\mathfrak g}^{\oplus 2}$
given by $y_{ij}\bigl((a_{pq}), (b_{rs})\bigr)=a_{ij}$,
$z_{ij}\bigl((a_{pq}), (b_{rs})\bigr)=b_{ij}$. Then
$y_{ij}=-y_{ji}$ and $z_{ij}=-z_{ji}$. Consider the
matrices $B:=(y_{ij})$, $C:=(z_{ij})$. Taking into
account that ${\rm tr}(PQ)={\rm tr}(QP)$ for any square
matrices $P$, $Q$, it is not difficult to deduce that
for every $\alpha_1, \alpha_2\in k$, the following
equalities hold:
\begin{equation}\label{coefficients1}
{\rm tr}((\alpha_1B+\alpha_2 C)^2)= \alpha_1^2{\rm
tr}(B^2)+2\alpha_1\alpha_2{\rm tr}(BC)+ \alpha_2{\rm
tr}(C^2),
\end{equation}
\begin{equation}\label{coefficients2}
\begin{split}
{\rm tr}((\alpha_1B+\alpha_2 C)^4)&= \alpha_1^4{\rm
tr}(B^4)+ 4\alpha_1^3\alpha_2{\rm tr}(B^3C)\\&+
2\alpha_1^2\alpha_2^2\bigl(2{\rm tr}(B^2C^2) + {\rm
tr}((BC)^2)\bigr)+ 4\alpha_1\alpha_2^3{\rm tr}(BC^3)+
\alpha_2^4{\rm tr}(C^4).
\end{split}
\end{equation}

 From \eqref{trace}, the definition of ${\rm
pol}_2k[\mathfrak g]^G$ (see Example
\ref{classical-setting1} and the first paragraph right
after~it), and \eqref{coefficients1},
\eqref{coefficients2}
 we deduce that ${\rm pol}_2k[\mathfrak g]^G$ is the
 algebra
$$k[{\rm tr}(B^2),
{\rm tr}(BC), {\rm tr}(C^2), {\rm tr}(B^4),
 {\rm tr}(B^3C), 2{\rm tr}(B^2C^2)+
{\rm tr}((BC)^2), {\rm tr}(BC^3), {\rm tr}(C^4)].$$
This shows that the transcendence degree of ${\rm
pol}_2k[\mathfrak g]^G$ over $k$ is not bigger than
$8$. On the other hand, since $\dim {\mathfrak
g}^{\oplus 2}=\dim G=20$, the transcendence degree of
$k[{\mathfrak g}^{\oplus 2}]^G$ over $k$ is not smaller
than $\dim {\mathfrak g}^{\oplus 2}-\dim G=10$, see,
e.g.,\;\cite[Theorem 3.3 and the Corollary of Lemma
2.4]{popov-vinberg} (actually it is equal to $10$
since, as one easily proves, the generic $G$-stabilizer
of the $G$-module ${\mathfrak g}^{\oplus 2}$ is
finite). Therefore $k[{\mathfrak g}^{\oplus 2}]^G$ is
not integral over ${\rm pol}_2k[\mathfrak g]^G$. By
\eqref{>0} and the Corollary~\ref{integral} we now
deduce that ${\rm pol\;ind}\;(\mathfrak g)=1$, i.e.,
$\mathfrak g$ is a $2$-algebra.

To complete the proof we have to calculate ${\rm
pol\hskip .7mm ind}\hskip .4mm(\mathfrak g)$ for
$\mathfrak g=\mathfrak g_1\oplus\ldots\oplus \mathfrak
g_m$ where $\mathfrak g_i=\mathfrak {sl}_2$ for every
$i$. We may assume that $G=G_1\times\ldots \times G_m$
where $G_i={\bf SL}_2$ for every $i$. The Corollary of
Lemma~\ref{weights} then reduces the proof to the case
$m=1$. Since the ${\bf SL}_2$-modules $\mathfrak
{sl}_2$ and $R_2$ are isomorphic, the claim now follows
from Theorem~\ref{sl2}. \quad $\square$
\renewcommand{\qed}{}
\end{proof}

Call a linear subspace $L$ of a reductive Lie algebra
$\mathfrak g$ {\it triangularizable} if there is
 a Borel subalgebra $\mathfrak b$ of
$\mathfrak g$ such that $L$ lies in the unipotent
radical of $\mathfrak b$ (for $\mathfrak g={\rm
Mat}_{n\times n}$, this means that $L$ is conjugate to
a subspace of the space of upper triangular matrices,
see\;\cite{gerstenhaber}, \cite{crt}, \cite{mor}). Call
$L$ {\it nilpotent} if every element of $L$ is
nilpotent.

\vskip 2mm

\begin{corollary} A semisimple Lie algebra $\mathfrak g$ contains a
$2$-dimensional nilpotent nontrian\-gu\-la\-rizable
linear subspace if and only if $\mathfrak g$ is not
isomorphic to $\mathfrak
{sl}_2\oplus\ldots\oplus\mathfrak {sl}_2$.
\end{corollary}


\providecommand{\bysame}{\leavevmode\hbox
to3em{\hrulefill}\thinspace}

 \end{document}